\documentclass[graybox]{svmult}

\usepackage{mathptmx}       
\usepackage{helvet}         
\usepackage{courier}        
\usepackage{type1cm}        
\usepackage{makeidx}         
\usepackage{graphicx}        
\usepackage{multicol}        
\usepackage[bottom]{footmisc}

\makeindex             

\usepackage{benstyle_book}
\usepackage{mathtools}
\usepackage{caption}
\usepackage{overpic}
\newcommand{\res}[1]{\ensuremath{#1\!\times\!#1}}

\DeclareGraphicsExtensions{.png, .jpg, .pdf}
\graphicspath{%
{./}%
{tmp/}%
{Diagrams/}%
{Diagrams/flip-rwt/}%
{Diagrams/flip-rwt2/}%
{Diagrams/brain/}%
{Diagrams/algs/}%
{Diagrams/changed/}%
{Diagrams/tvshuffle/}%
{Diagrams/tvshuffle2/}%
}

\begin{document}

\title*{The quest for optimal sampling: Computationally efficient, 
structure-exploiting measurements for compressed sensing}
\titlerunning{The quest for optimal sampling}
\author{Ben Adcock, Anders C.\ Hansen and Bogdan Roman}
\institute{Ben Adcock \at Purdue University, Department of Mathematics, 150 N.\ University St, West Lafayette, IN 47906, USA, \email{adcock@purdue.edu}
\and Anders C.\ Hansen \at DAMTP, Centre for Mathematical Sciences, University 
of Cambridge, Wilberforce Rd, Cambridge CB3 0WA, United Kingdom 
\email{a.hansen@maths.cam.ac.uk}
\and Bogdan Roman \at DAMTP, Centre for Mathematical Sciences, University of 
Cambridge, Wilberforce Rd, Cambridge CB3 0WA, United Kingdom 
\email{b.roman@maths.cam.ac.uk}
}
%
%

\maketitle



\abstract{
An intriguing phenomenon in many instances of compressed sensing is that the reconstruction quality is governed not just by the overall sparsity of the signal, but also on its structure.  This paper is about understanding this phenomenon, and demonstrating how it can be fruitfully exploited by the design of suitable sampling strategies in order to outperform more standard compressed sensing techniques based on random matrices.
}

\section{Introduction}
\label{s:introduction}

Compressed sensing concerns the recovery of signals and images from a small collection of linear measurements.  It is now a substantial area of research, accompanied by a mathematical theory that is rapidly reaching a mature state.  Applications of compressed sensing can roughly be divided into two areas.  First, \textit{type I} problems, where the physical device imposes a particular type of measurements.  This is the case in numerous real-world problems, including medical imaging (e.g.\ Magnetic Resonance Imaging (MRI) and Computerized Tomography (CT)), electron 
microscopy, seismic tomography and radar.  Second, \textit{type II} problems, where the sensing mechanism allows substantial freedom to design the measurements so as to improve the reconstructed image or signal.  Applications include compressive imaging and fluoresence microscopy.

This paper is devoted to the role of structured sparsity in both classes of problems.  It is well known that the standard sparsifying transforms of compressed sensing, i.e.\ wavelets and their various generalizations, not only give sparse coefficients, but that there is a distinct structure to this sparsity: wavelet coefficients of natural signals and images are far more sparse at fine scales than at coarse scales.  For type I problems, recent developments \cite{AHPRBreaking,BAGSAIEP} have shown that this structure plays a key role in the observed reconstruction quality.  Moreover, the optimal subsampling strategy depends crucially on this structure.  We recap this work in this paper.

Since structure is vitally important in type I problems, it is natural to ask whether or not it can be exploited to gain improvements in type II problems.  In this paper we answer this question in the affirmative.  We show that appropriately-designed structured sensing matrices can successfully exploit structure.  Doing so leads to substantial gains over classical approaches in type II problems, based on  convex optimization with universal, random Gaussian or Bernoulli measurements, as well as more recent structure-exploiting algorithms -- such as model-based compressed sensing \cite{BaranuikModelCS}, TurboAMP \cite{TurboAMP}, and Bayesian compressed sensing \cite{HeCarinStructCS,HeCarinTreeCS} -- which incorporate structure using bespoke recovery algorithms.  The proposed matrices, based on appropriately subsampled Fourier/Hadamard transforms, are also computationally efficient

We also review the theory of compressed sensing introduced in \cite{AHPRBreaking,Bastounis} with such structured sensing matrices.  The corresponding sensing matrices are highly non-universal and do not satisfy a meaningful Restricted Isometry Property (RIP).  Yet, recovery is still possible, and is vastly superior to that obtained from standard RIP matrices.  It transpires that the RIP is highly undesirable if one seeks to exploit structure by designing appropriate measurements.  Thus we consider an alternative that takes such structure into account, known as the \textit{RIP in levels} \cite{Bastounis}.

\section{Recovery of wavelet coefficients}
In many applications of compressed sensing, we are faced with the problem of recovering an image or signal $x$, considered as a vector in $\bbC^n$ or a function in $L^2(\mathbb{R}^d)$, that is sparse or compressible in an orthonormal basis of wavelets.  If $\Phi \in \bbC^{n \times n}$ or $\Phi \in \mathcal{B}(L^2(\mathbb{R}^d), \ell_2(\mathbb{N}))$ (the set of bounded linear operators) is the corresponding sparsifying transformation, then we write $x = \Phi c$, where $c \in \bbC^{n}$ or $c \in \ell_2(\mathbb{N})$ is the corresponding sparse or compressible vector of coefficients.  Given a sensing operator $A \in \bbC^{m \times n}$ or $A \in \mathcal{B}(L^2(\mathbb{R}^d), \mathbb{C}^m)$ and noisy measurements $y = A x + e$ with $\| e \|_2 \leq \eta$, the usual approach is to solve the $\ell_1$-minimization problem:
\be{
\label{general_CS}
\min_{z \in \bbC^n} \| \Phi z \|_1\quad \mbox{s.t.}\quad \| y - A z \|_2 \leq \eta.
}
or
\be{
\label{general_CS_Inf}
\inf_{z \in L^2(\mathbb{R}^d)} \| \Phi z \|_1\quad \mbox{s.t.}\quad \| y - A z \|_2 \leq \eta.
}
Throughout this paper, we shall denote a minimizer of \R{general_CS} or \R{general_CS_Inf} as $\hat{x}$.  Note that \R{general_CS_Inf} must be discretized in order to be solved numerically, and this can be done by restricting the minimization to be taken over a finite-dimensional space spanned by the first $n$ wavelets, where $n$ is taken sufficiently large \cite{BAACHGSCS}.

As mentioned, compressed sensing problems arising in applications can be divided into two classes:
\begin{enumerate}
\item[I.] \textit{Imposed sensing operators}.  The operator $A$ is specified 
by the practical device and is therefore considered fixed.  This is the 
case in MRI -- where $A$ arises by subsampling the Fourier transform 
\cite{Lustig3,Lustig}  -- as well as other examples, including X-ray CT 
(see \cite{ChartrandCSCTAlg} and references therein), radar 
\cite{HermanStrohmerRadar}, electron microscopy \cite{CSEIT,leary2013etcs}, seismic 
tomography 
\cite{LinHermanCSSeismic} and radio interferometry 
\cite{CSRadioInterferometry}. 
\item[II.] \textit{Designed sensing operators.}  The sensing mechanism 
allows substantial freedom to design $A$ so as to improve the compressed sensing 
reconstruction. Some applications belonging to this class are 
\textit{compressive imaging}, e.g.\ the single-pixel camera 
\cite{SinglePixelCamera} and the 
more recent lensless camera \cite{Lensless}, and compressive fluorescence 
microscopy \cite{Candes_PNAS}. In these applications $A$ is assumed to 
take binary values (typically $\{-1,1\}$), yet, as we will see later, this 
is not a significant practical restriction. 
\end{enumerate}

As stated, the purpose of this paper is to show that insight gained from understanding the application of compressed sensing to type I problems leads to more effective strategies for type II problems.

\subsection{Universal sensing matrices}
Let us consider type II problems.  In finite dimensions, the traditional 
compressed sensing approach has been to construct 
matrices $A$ possessing the following two properties.  First, they should 
satisfy the \textit{Restricted Isometry Property (RIP)}.  Second, they should be 
\textit{universal}.  That is, if $\Phi \in \bbC^{n \times n}$ is an arbitrary 
isometry, then $A \Phi$ also satisfies the RIP of the same order as $A$.  
Subject to these conditions, a typical result in compressed sensing is as 
follows (see  \cite{FoucartRauhutCSbook}, for example): if $A$ satisfies the 
RIP of order $2k$ with constant $\delta_{2k} < 4/\sqrt{41}$ then, for any $x 
\in \bbC^n$, we have 
\be{
\label{RIP_bound}
\| x - \hat{x} \|_2 \leq C \frac{\sigma_k(\Phi^* x)_1}{\sqrt{k}} + D \eta,
}
where $\hat{x}$ is any minimizer of \R{general_CS}, $C$ and $D$ are positive constants depending only on $\delta_{2k}$ and, for $c \in \bbC^n$, 
\bes{
\sigma_k(c)_1 = \inf_{z \in \Sigma_k} \| c - z \|_1,\qquad \Sigma_k = \{ z \in \bbC^n : \| z \|_0 \leq k \}.
}
Hence, $x$ is recovered exactly up to the noise level $\eta$ and the error $\sigma_k(c)_1$ of the best approximation of $c = \Phi^* x$ with a $k$-sparse vector.  Since $A$ is universal, one has complete freedom to choose the sparsifying transformation $\Phi$ so as to minimize the term $\sigma_k(\Phi^* x)_1$ for the particular signal $x$ under consideration.  


Typical examples of universal sensing matrices $A$ arise from random 
ensembles.  In particular, Gaussian or Bernoulli random matrices (with the 
latter having the advantage of being binary) both have this property with high 
probability whenever $m$ is proportional to $k$ times by a log factor.  For 
this reason, such matrices are often thought of as `optimal' matrices for compressed 
sensing.

\begin{remark}
One significant drawback of random ensembles, however, is that the 
corresponding matrices are dense and unstructured.  Storage and the lack of 
fast transforms render them impractical for all but small problem sizes.  To 
overcome this, various structured random matrices have also been developed and 
studied e.g.\ pseudo-random permutations of columns of Hadamard or DFT (Discrete Fourier Transform)
matrices \cite{HadamardBlockScrambled,Lensless}.  Often these admit fast, $\ord{n \log n}$ transforms.  
However, the 
best known theoretical RIP guarantees are usually larger than for 
(sub)Gaussian random matrices \cite{FoucartRauhutCSbook}.
\end{remark}

\subsection{Sparsity structure dependence and the flip test}

\begin{figure}[t]
\small
\begin{center}
\begin{tabular}{@{}c@{\hspace{0.005\linewidth}}c@{\hspace{0.005\linewidth}}c@{}}
\includegraphics[width=0.33\linewidth]{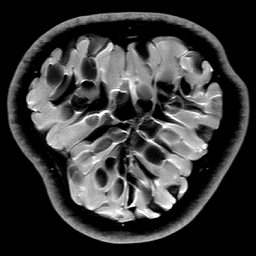}&
\includegraphics[width=0.33\linewidth]{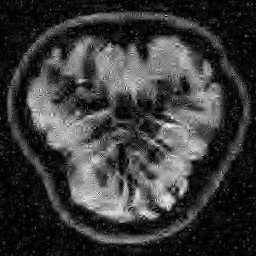}&
\includegraphics[width=0.33\linewidth]{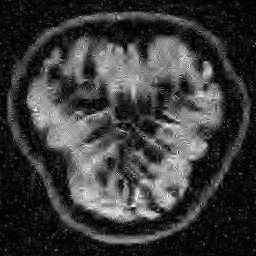}\\
Original image & Gauss. to DB4, err=31.54\% & Gauss. to Flip DB4, err=31.51\% 
\\[10pt]
\includegraphics[width=0.33\linewidth]{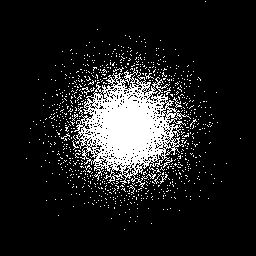}&
\includegraphics[width=0.33\linewidth]{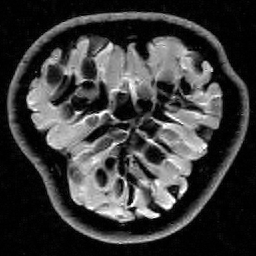}&
\includegraphics[width=0.33\linewidth]{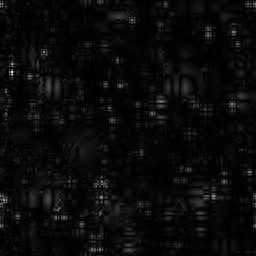}\\
Subsampling pattern & DFT to DB4, err=10.96\% & DFT to Flipped DB4, 
err=99.3\% \\
\end{tabular}
\caption{Recovery of an MRI image of a passion fruit using \R{general_CS} from 
$m =8192$ samples at $n=256 \times 256$ resolution (i.e.\ a 12.5\% subsampling 
rate) with Daubechies-4 wavelets as the sparsifying transform $\Phi$. Top row: 
Flip test for Gaussian random measurements.  Botton row: Flip test for 
subsampled DFT measurements taken according to the subsampling pattern shown 
in the bottom left panel.  The flip test suggests that the RIP holds for 
random sensing matrices (top row), but that there is no RIP for structured 
sensing matrices with structured sampling (bottom row).}
\label{f:BernStoneHad}
\end{center}
\end{figure}

Since it will become important later, we now describe a quick and simple test, 
which we call the \textit{flip test}, to investigate the 
presence or absence of an RIP.  Success of this test suggests the existence of 
an RIP and failure demonstrates its lack.

Let $A \in \bbC^{m \times n}$ be a sensing matrix,  $x \in \bbC^{n}$ an image and $\Phi \in \bbC^{n \times n}$ a sparsifying transformation.  Recall that sparsity of the vector $c = \Phi^* x$ is unaffected by permutations.  Thus, let us define the flipped vector
\bes{
P(c) = c' \in \bbC^n,\quad c'_i = c_{n+1-i},\quad i=1,\ldots,n,
}
and using this, we construct the flipped image
$
x' = \Phi c'.
$
Note that, by construction, we have $\sigma_k(c)_1 = \sigma_k(c')_1$.  Now suppose we perform the usual compressed sensing reconstruction \R{general_CS} on both $x$ and $x'$, giving approximations $\hat{x} \approx x$ and $\hat{x}' \approx x'$.  We now wish to reverse the flipping operation.  Thus, we compute
$
\check{x} = \Phi P(\Phi^* \hat{x}'),
$
which gives a second approximation to the original image $x$.  

This test provides a simple way to investigate whether or not the RIP holds.  To see why, suppose that $A$ satisfies the RIP.  Then by construction, we have that
\bes{
\| x - \hat{x} \|_2 , \| x - \check{x} \|_2 \leq C \frac{\sigma_k(\Phi^* x)_1}{\sqrt{k}} + D \eta.
}
Hence both $\hat{x}$ and $\check{x}$ should recover $x$ equally well.  In the top row of Figure \ref{f:BernStoneHad} we present the result of the flip test for a Gaussian random matrix.  As is evident, the reconstructions $\hat{x}$ and $\check{x}$ are comparable, thus indicating the RIP.  

Having considered type II problems, let us now examine the flip test for a type I problem.  As discussed, in applications such as MRI, X-ray CT, radio interferometry, etc, the matrix $A$ is imposed by the physical sensing device and arises from subsampling the rows of the DFT matrix $F \in \bbC^{n \times n}$.\footnote{In actual fact, the sensing device takes measurements of the \textit{continuous} Fourier transform of a function $x \in L^2(\bbR^d)$.  As discussed in \cite{BAACHGSCS,BAGSAIEP}, modelling continuous Fourier measurements as discrete Fourier measurements can lead to inferior reconstructions, and worse, inverse crimes.  To avoid this, one must consider an infinite-dimensional compressed sensing approach, as in \R{general_CS_Inf}.  See \cite{AHPRBreaking,BAGSAIEP} for details, as well as \cite{PruessmannUnserMRIFast} for implementation in MRI.  However, for simplicity, we shall continue to work with the finite-dimensional model in the remainder of this paper.}
  Whilst one often has some freedom to choose which rows to sample (corresponding to selecting particular frequencies at which to take measurements), one cannot change the matrix $F$.

It is well known that in order to ensure a good reconstruction, one cannot subsample the DFT uniformly at random (recall that the sparsifying transform is a wavelet basis), but rather one must sample randomly according to an appropriate nonuniform density \cite{AHPRBreaking,Candes_Romberg,Lustig,WangAcre}.  See the bottom left panel of Figure \ref{f:BernStoneHad} for an example of a typical density.  As can be seen in the next panel, by doing so one achieves a great recovery.  However, the result of the flip test in the bottom right panel clearly demonstrates that the matrix $F \Phi$ does not satisfy an RIP.  In particular, the ordering of the wavelet coefficients plays a crucial role in the reconstruction quality.  To explain this, and in particular, the high-quality reconstruction seen in the unflipped case, one evidently requires a new analytical framework. 

Note that the flip test in Figure \ref{f:BernStoneHad} also highlights another 
important phenomenon: namely, the effectiveness of the subsampling strategy 
depends on the sparsity structure of the image. In particular, two images with 
the same total sparsity  (the original $x$ and the flipped $x'$) result in wildly different errors when the same 
sampling pattern is used.  Thus we conclude that there is no one optimal sampling strategy for all sparse vectors of wavelet coefficients.

\begin{figure}[t]
\small
\begin{center}
\begin{tabular}{@{}c@{\hspace{0.005\linewidth}}c@{\hspace{0.005\linewidth}}%
c@{\hspace{0.005\linewidth}}c@{\hspace{0.005\linewidth}}c@{}}
\includegraphics[width=0.195\linewidth]{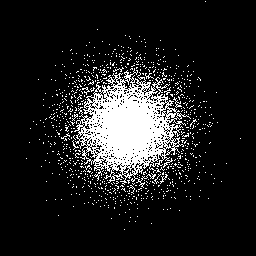}&
\includegraphics[width=0.195\linewidth]{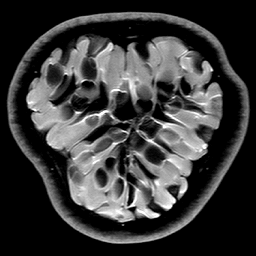}&
\includegraphics[width=0.195\linewidth]{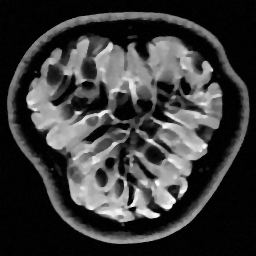}&
\includegraphics[width=0.195\linewidth]{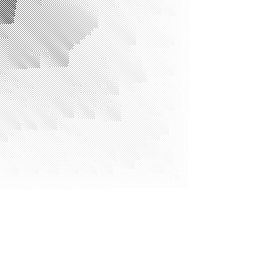}&
\includegraphics[width=0.195\linewidth]{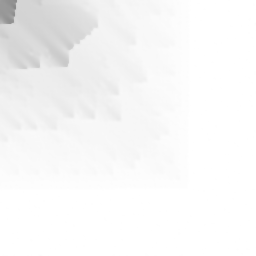}\\
Subsample pattern & Original & TV, \textbf{err = 9.27\%} & Permuted gradients 
& TV, \textbf{err = 3.83\%}\\[5pt]
\includegraphics[width=0.195\linewidth]{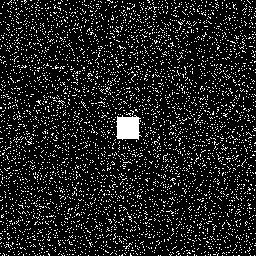}&
\includegraphics[width=0.195\linewidth]{0256_12p5_dft_tv_0p02_ml_Z8L8_ref.png}&
\includegraphics[width=0.195\linewidth]{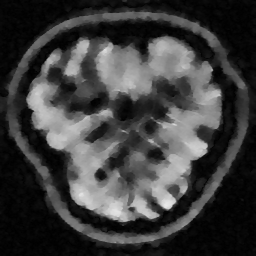}&
\includegraphics[width=0.195\linewidth]{0256_12p5_dft_tv_0p02_ml_equiv_DLZJ_ref.png}&
\includegraphics[width=0.195\linewidth]{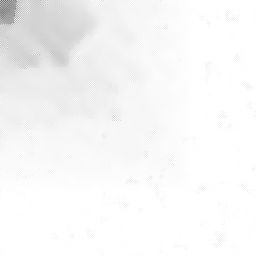}\\
Subsample pattern & Original & TV, \textbf{err = 28.78\%} & Permuted gradients 
& 
TV, \textbf{err = 3.76}\%
\end{tabular}
\caption{TV recovery \R{general_CS_TV} from $m = 16384$ DFT samples at $n=512 
\times 512$ (6.25\% subsampling rate).  The \textit{Permuted gradients} image 
was built from the same image gradient vectors as the \textit{Original} image, 
having the same TV norm and gradient sparsity structure, differing only in the 
ordering and sign of the gradient vectors. The large error differences confirm 
that, much like the flip test for wavelet coefficients, the sparsity structure 
matters for TV reconstructions as well.}
\label{f:TV}
\end{center}
\end{figure}

Let us also note that the same conclusions of 
the flip test hold when \R{general_CS} with wavelets is replaced by TV-norm 
minimization:
\be{
\label{general_CS_TV}
\min_{z \in \bbC^n} \| z \|_{TV} \quad \mbox{s.t.}\quad \| y - A z \|_2 \leq \eta.
}
Recall that $\|x\|_{TV} = \sum_{i,j} \|\nabla x(i,j)\|_2$, where we have
$
\nabla x(i,j) = 
\{D_1x(i,j), D_2x(i,j)\},$ $D_1x(i,j) = x(i+1,j)-x(i,j),$ $D_2x(i,j) = 
x(i,j+1)-x(i,j)$.  In the experiment leading to Figure \ref{f:TV}, we chose an image $x\in[0,1]^{N\times N}$, and then built a different image $x'$ from the gradient of $x$ so that $\{\|\nabla x'(i,j)\|_2\}$ is a permutation of $\{\|\nabla x(i,j)\|_2\}$ for which $x'\in[0,1]^{N\times N}$.  Thus, the two images have the same ``TV sparsity" and the same TV norm. In 
Figure \ref{f:TV} we demonstrate how the errors differ substantially for the 
two images when using the same sampling pattern. Note also how the improvement 
depends both on the TV sparsity structure and on the subsampling pattern. 
Analysis of this phenomenon is work in progress.  In the remainder of this paper we will focus 
on structured sparsity for wavelets.

\subsection{Structured sparsity}

\begin{figure}[t]
\begin{center}
\def\imagetop#1{\vtop{\null\hbox{#1}}}
\begin{tabular}{@{\hspace{0pt}}c@{\hspace{0.1\textwidth}}c@{\hspace{0pt}}}
\imagetop{\includegraphics[width=0.27\linewidth]{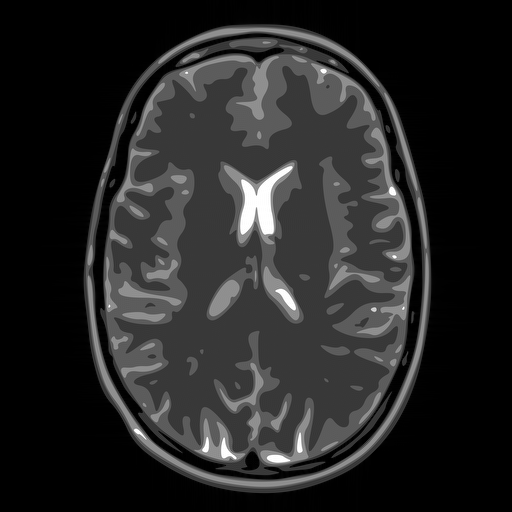}}&
\imagetop{\includegraphics[width=0.42\textwidth]{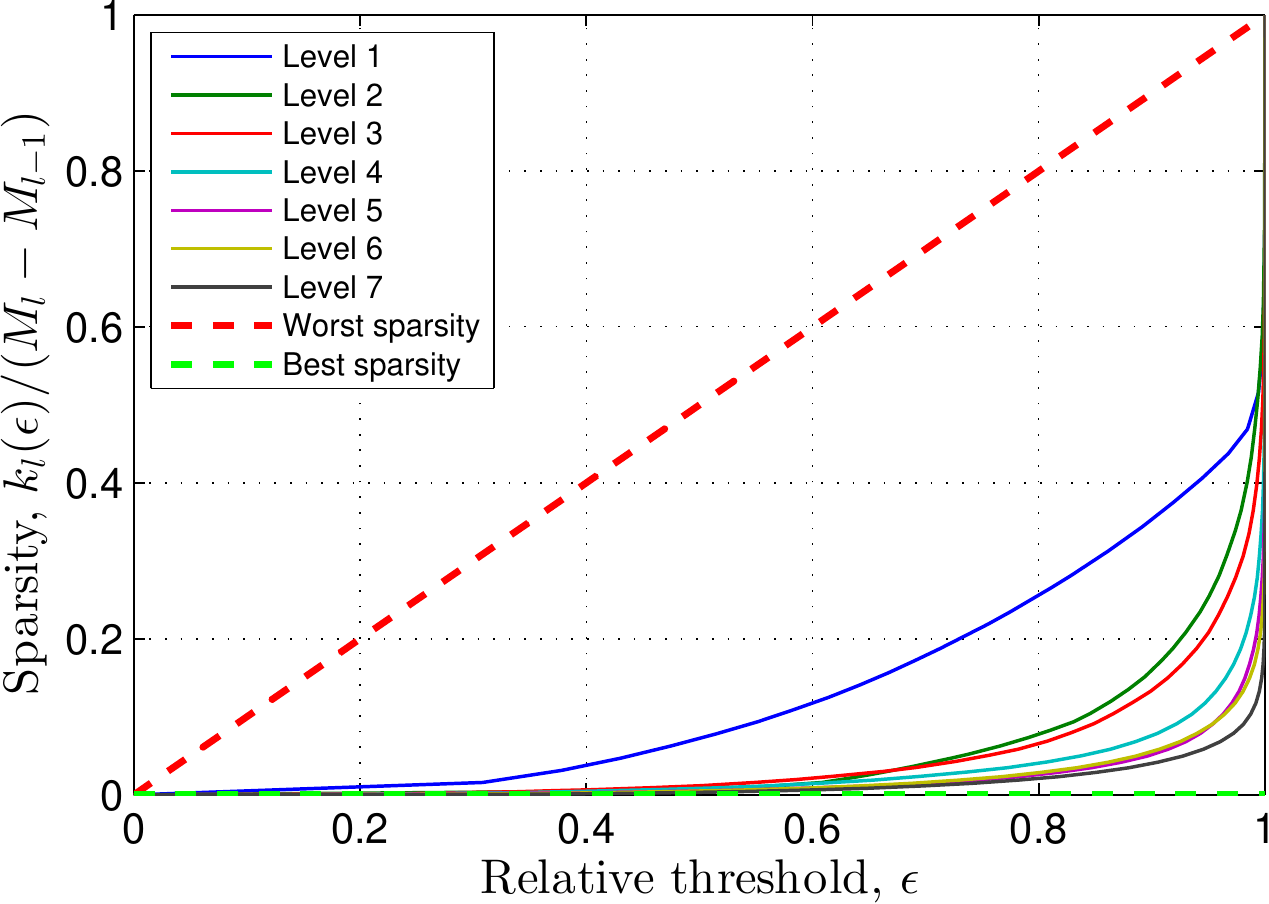}}
\end{tabular}
\caption{Left: GLPU phantom \cite{GLPU}.  Right: Relative sparsity 
of Daubechies-4 coefficients.  Here the levels correspond to wavelet scales 
and 
$k_l(\epsilon)$ is given by \R{sk_def1}.  Each curve shows the relative 
sparsity at level $l$ as a function of $\epsilon$.  The decreasing nature of 
the curves for increasing $l$ confirms \R{fine_decay}.
}
\label{fig:CS_LevelsSparsity}
\end{center}
\end{figure}

One of the foundational results of nonlinear approximation is that, for natural images and signals, the best $k$-term approximation error in a wavelet basis decays rapidly in $k$ \cite{DeVoreNLACTA,mallat09wavelet}.  In other words, wavelet coefficients are approximately $k$-sparse.  
However, wavelet coefficients possess far more structure than mere sparsity.  Recall that a wavelet basis for $L^2(\bbR^d)$ is naturally partitioned into dyadic scales.  Let $0 = M_0 < M_1 <\ldots < \infty$ be such a partition, and note that $M_{l+1} - M_{l} = \ord{2^l}$ in one dimension and $M_{l+1} - M_{l} = \ord{4^l}$ in two dimensions.  If $x = \Phi c$, let $c^{(l)} \in \bbC^{M_{l}-M_{l-1}}$ denote the wavelet coefficients of $x$ at scale $l=1,2,\ldots$, so that $c = ( c^{(1)} |  c^{(2)} | \cdots )^{\top}$.  Suppose that $\epsilon > 0$ and define
\be{
\label{sk_def1}
k_l = k_l(\epsilon) = \min \left \{ K : \sum^{K}_{i=1} |c^{(l)}_{\pi(i)} |^2 \geq \epsilon^2 \| c^{(l)} \|^2_2 \right \},\quad l=1,2,\ldots,
} 
where $\pi$ is a bijection that gives a nonincreasing rearrangement of the entries of $c^{(l)}$, i.e.\ $|c^{(l)}_{\pi(i)}| \geq |c^{(l)}_{\pi(i+1)} |$ for $i=1,\ldots,M_{l}-M_{l-1} - 1$.  Sparsity of the whole vector $c$ means that for large $r$ we have $k / n \ll 1$, where
$
k = k(\epsilon) = \sum^{r}_{l=1} k_l,
$
is the total effective sparsity up to finest scale $r$.  However, Figure \ref{fig:CS_LevelsSparsity} reveals that not only is this the case in practice, but we also have so-called \textit{asymptotic sparsity}.  That is
\be{
\label{fine_decay}
k_l / (M_{l}-M_{l-1}) \rightarrow 0,\quad l \rightarrow \infty.
}
Put simply, wavelet coefficients are much more sparse at fine scales than they are at coarse scales.

Note that this observation is by no means new: it is a simple consequence of 
the dyadic scaling of wavelets, and is a crucial step towards 
establishing the nonlinear approximation result mentioned above.  However, 
given that wavelet coefficients always exhibit such structure, one may ask the 
following question.  For type II problems, are the traditional sensing matrices 
of compressed sensing -- which, as shown by the flip test (top row of Figure 
\ref{f:BernStoneHad}), recover all sparse vectors of coefficients equally 
well, regardless of ordering -- optimal when wavelets are used as the 
sparsifying transform?  It has been demonstrated in Figure 
\ref{f:BernStoneHad} (bottom row) that structure plays a key role in type I compressed sensing problems.  Leveraging this insight, in the next section we show that significant gains are possible for type II problems in terms of both 
the reconstruction quality and computational efficiency when the sensing 
matrix $A$ is designed specifically to take advantage of such  inherent 
structure.

\begin{remark}
Asymptotic sparsity is by no means limited to wavelets.  A similar property holds for other types of -lets, including curvelets \cite{Cand,candes2004new}, contourlets \cite{Vetterli,Do} or shearlets \cite{Gitta2,Gitta,Gitta3} (see \cite{BAGSAIEP} for examples of Figure \ref{fig:CS_LevelsSparsity} based on these transforms).  More generally, \textit{any} sparsifying transform that arises (possibly via discretization) from a countable basis or frame will typically exhibit asymptotic sparsity.
\end{remark}

\begin{remark}\label{r:tree_structure}
Wavelet coefficients (and their various generalizations)  actually possess far 
more structure than the asymptotic sparsity \R{fine_decay}.  Specifically, 
they tend to live on rooted, connected trees \cite{CrouseEtAlWaveleTree}. There are a 
number of existing algorithms which seek to 
exploit such structure within a compressed sensing context.  We shall discuss 
these further in Section \ref{ss:Structured_Sampling_Recovery}.
\end{remark}

\section{Efficient sensing matrices for structured sparsity}
\label{s:SensMatSS}

For simplicity, we now consider the finite-dimensional setting, although the arguments extend to the infinite-dimensional case \cite{AHPRBreaking}.  Suppose that $\Phi \in \bbC^{n \times n}$ corresponds to a wavelet basis so that $c = \Phi^* x$ is not just $k$-sparse, but also asymptotically sparse with the sparsities $k_1,\ldots,k_r$ within the wavelet scales being known.  As before, write $c^{(l)}$ for set of coefficients at scale $l$.  Considering type II problems, we now seek a sensing matrix $A$ with as few rows $m$ as possible that exploits this local sparsity information.

\subsection{Block-diagonality and structured Fourier/Hadamard sampling}
For the $l^{\rth}$ level, suppose that we assign a total of $m_l \in \bbN$ rows of $A$ in order to recover the $k_l$ nonzero coefficients of $c^{(l)}$.  Note that $m = \sum^{r}_{l=1} m_l$.  Consider the product $B = A \Phi$ of the sensing matrix $A$ and the sparsifying transform $\Phi$.  Then there is a natural decomposition of $B$ into blocks $\{B_{jl}\}^{r}_{j,l=1}$ of size $m_j \times (M_l-M_{l-1})$, where each block corresponds to the $m_j$ measurements of the $M_l-M_{l-1}$ wavelet functions at the $l^{\rth}$ scale.

Suppose it were possible to construct a sensing matrix $A$ such that (i) $B$ was block diagonal, i.e.\ $B_{jl} = 0$ for $j \neq l$, and (ii) the diagonal blocks $B_{ll}$ satisfied an RIP of order $2 k_l$ whenever $m_l$ was proportional to $k_l$ times by the usual log factor.  In this case, one recovers the coefficients $c^{(l)}$ at scale $l$ from near-optimal numbers of measurements using the usual reconstruction \R{general_CS}.

This approach, originally proposed by Donoho \cite{donohoCS} and Tsaig \& 
Donoho \cite{DonohoCSext} under the name of `multiscale compressed sensing', 
allows for structure to be exploited within compressed sensing.  Similar ideas 
were also pursued by Romberg \cite{RombergCompImg} within the context of 
compressive imaging.  Unfortunately, it is normally impossible to design an $m 
\times N$ matrix $A$ such that $B = A \Phi$ is exactly block diagonal.  
Nevertheless, the notion of block-diagonality provides insight into better 
designs for $A$ than purely random ensembles.  To proceed, we relax the 
requirement of strict block-diagonality, and instead ask whether there exist 
practical sensing matrices $A$ for which $B$ is approximately block-diagonal 
whenever the sparsifying transform $\Phi$ corresponds to wavelets.  
Fortunately, the answer to this question is affirmative: as we shall explain 
next, and later confirm with a theorem, approximate block-diagonality can be 
ensured whenever $A$ arises by appropriately subsampling the rows of the Fourier or Hadamard transform.  Recalling that the former arises naturally in type I problems, this points towards the previously-claimed conclusion that new insight brought 
about by studying imposed sensing matrices leads to better approaches for the 
type II problem.

Let $F \in \bbC^{n \times n}$ be either the discrete Fourier or discrete Hadamard transform.  Let $\Omega \subseteq \{1,\ldots,n\}$ be an index set of size $| \Omega | =m$.  We now consider choices for $A$ of the form
$
A = P_{\Omega} F,
$
where $P_{\Omega} \in \bbC^{m \times n}$ is the restriction operator that selects rows of $F$ whose indices lie in $\Omega$.  We seek an $\Omega$ that gives the desired block-diagonality.  To do this, it is natural to divide up $\Omega$ itself into $r$ disjoint blocks
\bes{
\Omega = \Omega_1 \cup \cdots \cup \Omega_r,\quad | \Omega _l| = m_l,
}
where the $l^{\rth}$ block $\Omega_l \subseteq \{ N_{l-1},\ldots,N_l \}$ corresponds to the $m_l$ samples required to recover the $k_l$ nonzero coefficients at scale $l$.  Here the parameters $0=N_0 < N_1 < \ldots < N_r = n$ are appropriately chosen and delineate frequency bands from which the $m_l$ samples are taken.

In Section \ref{ss:Fourier_localinc}, we explain why this choice of $A$ works, and in particular, how to choose the sampling blocks $\Omega_l$.  In order to do this, it is first necessary to recall the notion of incoherent bases.

\subsection{Incoherent bases and compressed sensing}\label{ss:inc_CS}

Besides random ensembles, a common approach in compressed sensing is to design sensing matrices using orthonormal systems that are incoherent with the particular choice of sparsifying basis $\Phi$ \cite{Candes_Plan,Candes_Romberg,FoucartRauhutCSbook}.  Let $\Psi \in \bbC^{n \times n}$ be an orthonormal basis of $\bbC^n$.  The (mutual) coherence of $\Phi$ and $\Psi$ is the quantity
\bes{
\mu = \mu(\Psi^*\Phi) = \max_{i,j=1,\ldots,n} | (\Psi^* \Phi)_{i,j} |^2.
}
We say $\Psi$ and $\Phi$ are \textit{incoherent} if $\mu(\Psi,\Phi) \leq a/ n$ for some $a\geq1$ independent of $n$.  Given such a $\Psi$, one constructs the sensing matrix
$
A = P_{\Omega} \Psi,
$
where $\Omega \subseteq \{1,\ldots,N \}$, $| \Omega | = m$ is chosen uniformly at random.  A standard result gives that a $k$-sparse signal $x$ in the basis $\Phi$ is recovered exactly with probability at least $1-p$, provided
\bes{
m \geq C k a \log(1+p^{-1}) \log(n),
}
for some universal constant $C >0$ \cite{FoucartRauhutCSbook}.  As an example, consider the Fourier basis $\Psi = F$.  This is incoherent with the canonical basis $\Phi = I$ with optimally-small constant $a=1$.  Fourier matrices subsampled uniformly at random are efficient sensing matrices for signals that are themselves sparse.

However, the Fourier matrix is not incoherent with a wavelet basis: $\mu(F , 
\Phi) = \ord{1}$ as $n \rightarrow \infty$ for any orthonormal wavelet basis 
\cite{AHPRBreaking}.  Nevertheless, Fourier samples taken within appropriate 
frequency bands (i.e.\ not uniformly at random) are \textit{locally} 
incoherent with wavelets in the corresponding scales.  This observation, which 
we demonstrate next, explains the success of the sensing matrix $A$ 
constructed in the previous subsection for an appropriate choice of 
$\Omega_1,\ldots,\Omega_r$.

\subsection{Local incoherence and near block-diagonality of Fourier measurements with wavelets}\label{ss:Fourier_localinc}
For expository purposes, we consider the case of one-dimensional Haar wavelets.  We note however that the arguments generalize to arbitrary compactly-supported orthonormal wavelets, and to the infinite-dimensional setting where the unknown image $x$ is a function.  See Section \ref{ss:Fourier_wavelet}.  

Let $j=0,\ldots,r-1$ (for convenience we now index from $0$ to $r-1$, as 
opposed to $1$ to $r$) be the scale and $p=0,\ldots,2^{j}-1$ the translation.  The Haar basis consists of the functions $\{ \psi \} \cup \{ \phi_{j,p} 
\}$, where $\psi \in \bbC^n$ is the normalized scaling function and 
$\phi_{j,p}$ are the scaled and shifted versions of the mother wavelet $\phi 
\in \bbC^n$.  It is a straightforward, albeit tedious, exercise to show that
\bes{
| \cF\phi_{l,p}(\omega) | = 2^{l/2-r+1} \frac{| \sin(\pi \omega / 2^{l+1}) |^2}{| \sin (\pi \omega / 2^r ) |} \lesssim 2^{l/2} \frac{| \sin(\pi \omega / 2^{l+1}) |^2}{| \omega |},\quad | \omega | < 2^r,
}
where $\cF$ denotes the DFT \cite{AHRHaarFourierNote}.  This suggests that the Fourier transform $\cF \phi_{l,p}(\omega)$ is large when $\omega \approx 2^l$, yet smaller when $\omega \approx 2^j$ with $j \neq l$.  Hence we should separate frequency space into bands of size roughly $2^j$. 

Let $F \in \bbC^{n \times n}$ be the DFT matrix with rows indexed from $-n/2+1$ to $n/2$.  Following an approach of \cite{Candes_Romberg}, we now divide these rows into the following disjoint frequency bands
\bes{
W_0 = \{ 0,1 \},\quad W_{j} = \{ - 2^{j}+1,\ldots,-2^{j-1} \} \cup \{ 2^{j-1}+1,\ldots,2^{j} \},\quad l=0,\ldots,{r-1}.
}
With this to hand, we now define $\Omega_{j}$ to be a subset of $W_j$ of size $| \Omega_j | = m_j$ chosen uniformly at random.  Thus, the overall sensing matrix $A = P_{\Omega} F$ takes measurements of the signal $x$ by randomly drawing $m_j$ samples its Fourier transform within the frequency bands $W_j$.

Having specified $\Omega$, let us note that the matrix $U = F \Phi$ naturally divides into blocks, with the rows corresponding to the frequency bands $W_l$ and the columns corresponding to the wavelet scales.  Write $U = \{ U_{jl} \}^{r-1}_{j,l=0}$ where $U_{jl} \in \bbC^{2^{j} \times 2^{l}}$.    For compressed sensing to succeed in this setting, we require two properties.  First, the diagonal blocks $U_{jj}$ should be incoherent, i.e.\ $\mu(U_{jj}) = \ord{2^{-j}}$.  Second, the coherences $\mu(U_{jl})$ of the off-diagonal blocks $U_{jl}$ should be appropriately small in comparison to $\mu(U_{jj})$.  These two properties are demonstrated in the following lemma:

\lem{
\label{l:DFT_Haar_local_coherence}
We have $\mu(U_{jj}) \lesssim 2^{-j}$ and, in general,
\bes{
\mu(U_{jl}) \lesssim \mu(U_{jj}) 2^{-|j-l|},\quad j,l=0,\ldots,r-1
}
}
Hence $U$ is approximately block diagonal, with exponential decay away from the diagonal blocks.  Fourier measurements subsampled according to the above strategy are therefore ideally suited to recover structured sparse wavelet coefficients.\footnote{For brevity, we do not give the proof of this lemma or the later recovery results for Haar wavelets, Theorem \ref{t:Fourier_Haar}.  Details of the proof can be found in the short note \cite{AHRHaarFourierNote}.}

\begin{figure}[t]
\centering
\includegraphics[width=0.325\textwidth]{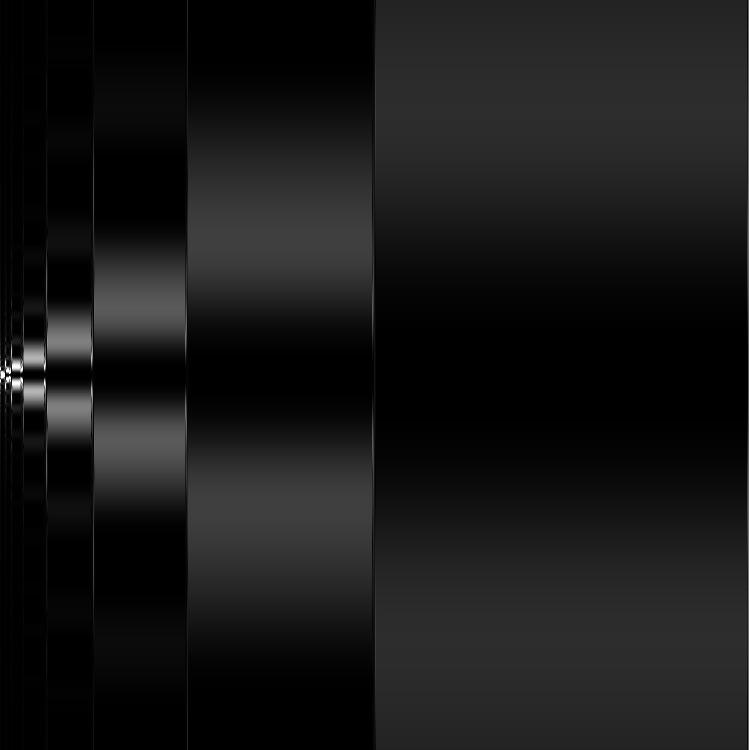}
\includegraphics[width=0.325\textwidth]{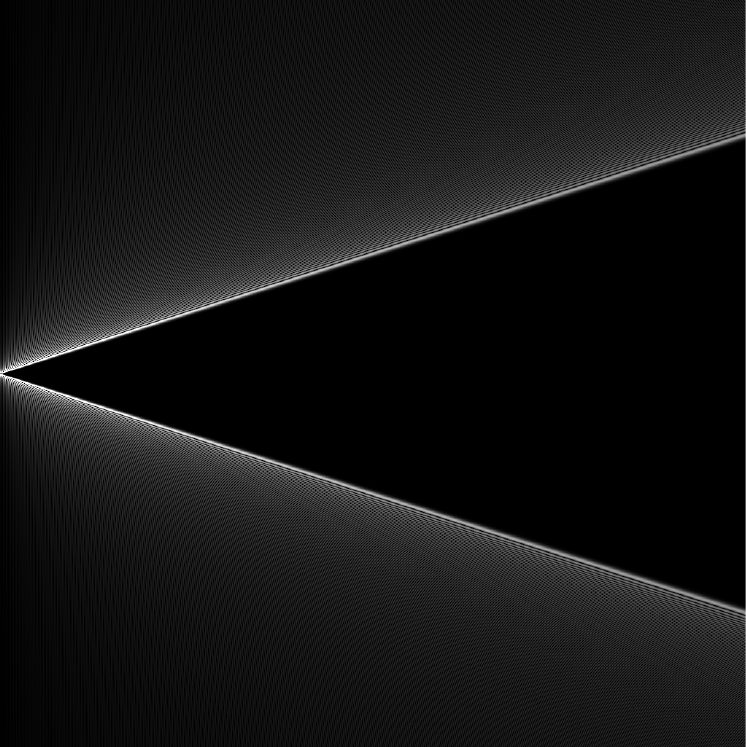}
\includegraphics[width=0.325\textwidth]{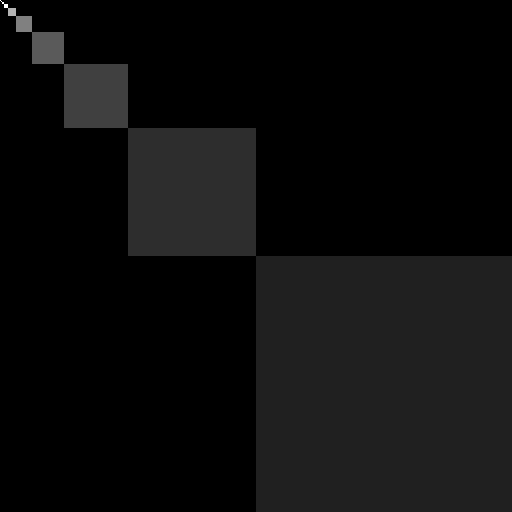}
\caption{The absolute values of the matrix $U = F \Phi$, where $F$ is the 
discrete Fourier transform (left and middle) or the Hadamard transform (right) 
and the sparsifying transform $\Phi$ corresponds to Haar wavelets (left and right) or Legendre polynomials (middle).  Light colours correspond to large 
values and dark colours to small values.}
\label{fig:U_plot}       
\end{figure}

The left panel of Figure \ref{fig:U_plot} exhibits this decay by plotting the 
absolute values of the matrix $U$.  In the right panel, we also show a similar 
result when the Fourier transform is replaced by the Hadamard transform. This is 
an important case, since the measurement matrix is binary.  The middle panel 
of the figure shows the $U$ matrix when Legendre polynomials are used as the 
sparsifying transform, as is sometimes the case for smooth signals.  It 
demonstrates that diagonally-dominated coherence is not just a phenomenon 
associated with wavelets.

\begin{figure}
\centering
\begin{tabular}{@{}c@{\hspace{0.00875\textwidth}}c@{\hspace{0.00875\textwidth}}c@{}}
\includegraphics[width=0.325\textwidth]{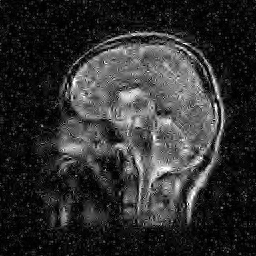}&
\includegraphics[width=0.325\textwidth]{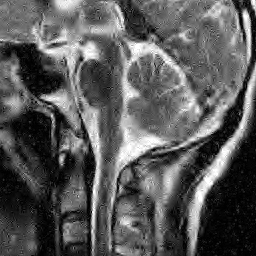}&
\includegraphics[width=0.325\textwidth]{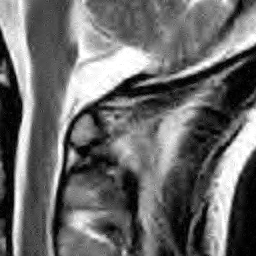}\\
\res{256}, err = 41.6\% & \res{512}, err = 25.3\% & \res{1024}, err = 11.6\% 
\\[5pt]
\includegraphics[width=0.325\textwidth]{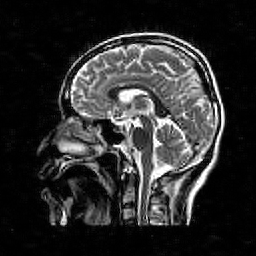}&
\includegraphics[width=0.325\textwidth]{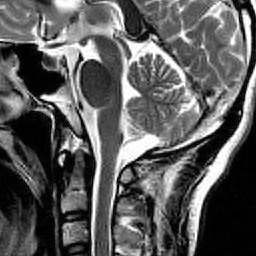}&
\includegraphics[width=0.325\textwidth]{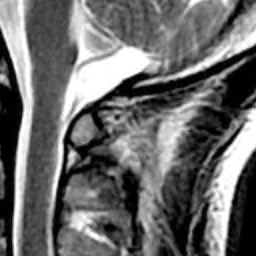}\\
\res{256}, err = 21.9\% & \res{512}, err = 10.9\% & \res{1024}, err = 3.1\%
\end{tabular}
\caption{Recovery from 12.5\% measurements using \R{general_CS} with 
Daubechies-4 wavelets as the sparsifying transform.  Top row: Random Bernoulli 
sensing matrix. Bottom row: Fourier sensing matrix with multilevel subsampling (see 
Definition \ref{multi_level_dfn}). All images are \res{256} crops of the 
original full resolution versions in order to aid the visual comparison.}
\label{fig:GaussMultilevelComp}       
\end{figure}

Having identified a measurement matrix to exploit structured sparsity, let us demonstrate its effectiveness.  In Figure 
\ref{fig:GaussMultilevelComp} we compare these measurements with the case of 
random Bernoulli measurements (this choice was made over random Gaussian 
measurements because of storage issues).  As is evident, at all resolutions we 
see a significant advantage, since the former strategy exploits the structured 
sparsity.  Note that for both approaches, the reconstruction quality is \textit{resolution dependent}: the error decreases as the resolution increases, due to the increasing sparsity of wavelet coefficients at higher resolutions.  However, because the Fourier/wavelets matrix $U$ is \textit{asymptotically incoherent} (see also Section \ref{ss:concepts}), it exploits the inherent asymptotic sparsity structure \R{fine_decay} of the wavelet coefficients as the resolution increases, and thus gives successively greater improvements over random Bernoulli measurements.

\begin{remark}\label{r:comp_time}
Besides improved reconstructions, an important feature of this approach is 
storage and computational time.  Since $F$ and $\Phi$ have fast, $\ord{n \log 
n}$, transforms, the matrix $A = P_{\Omega} F \Phi$ does not need to be 
stored, and the reconstruction \R{general_CS} can be performed efficiently 
with standard $\ell_1$ solvers (we use SPGL1 \cite{SPGL1,SPGL1Paper} 
throughout). 
\end{remark}

Recall that in type I problems such as MRI, we are constrained by the physics 
of the device to take Fourier measurements.  A rather strange conclusion of 
Figure \ref{fig:GaussMultilevelComp} is the following: compressed sensing 
actually works better for MRI with the intrinsic measurements, than if one 
were able to take optimal (in the sense of the standard sparsity-based theory) 
random (sub)Gaussian measurements.  This has practical consequences.  In MRI 
there is actually a little flexibility to design measurements, based on specifying
appropriate pulses.  By doing this, a number of approaches
\cite{HaldarRandomEncoding,ToeplitzMRI,VanderEtAlSpreadSpectrum,SebertRandomMRI,WongMRI} have been proposed 
to make MRI measurements closer to uniformly incoherent with wavelets (i.e.\ similar to 
random Gaussians).  On the other hand, Figure \ref{fig:GaussMultilevelComp} suggests that one can obtain great results in practice by appropriately subampling the unmodified Fourier operator.

\section{A general framework for compressed sensing based on structured sparsity}

Having argued for the particular case of Fourier samples with Haar wavelets, we now describe a general mathematical framework for structured sparsity.  This is based on work in \cite{AHPRBreaking}.

\subsection{Concepts}\label{ss:concepts}
We shall work in both the finite- and infinite-dimensional settings, where $U \in \bbC^{n \times n}$ or $U \in \cB(\ell_2(\bbN))$ respectively.  We assume throughout that $U$ is an isometry.  This occurs for example when $U = \Psi^* \Phi$ for an orthonormal basis $\Psi$ and an orthonormal system $\Phi$, as is the case for the example studied in the previous section: namely, Fourier sampling and a wavelet sparsifying basis.  However, framework we present now is valid for an arbitrary isometry $U$, not just this particular example.  We discuss this case further in Section \ref{ss:Fourier_wavelet}.

We first require the following definitions.  In the previous section it was suggested to divide both the sampling strategy and the sparse vector of coefficients into disjoint blocks.  We now formalize these notions:

\defn{[Sparsity in levels]
\label{d:Asy_Sparse}
Let $c$ be an element of either $\bbC^N$ or $\ell_2(\bbN)$. For $r \in \bbN$ let $\mathbf{M} = (M_1,\ldots,M_r) \in \bbN^r$ with $1 \leq M_1 < \ldots < M_r$ and $\mathbf{k} = (k_1,\ldots,k_r) \in \bbN^r$, with $k_j \leq M_j - M_{j-1}$, $k=1,\ldots,r$, where $M_0 = 0$.  We say that $c$ is $(\mathbf{k},\mathbf{M})$-sparse if, for each $j=1,\ldots,r$, the set
\bes{
\Delta_j  : = \mathrm{supp}(c) \cap \{ M_{j-1}+1,\ldots,M_{j} \},
}
satisfies $| \Delta_j | \leq k_j$.  We denote the set of $(\mathbf{k},\mathbf{M})$-sparse vectors by $\Sigma_{\mathbf{k},\mathbf{M}}$.
}

This definition allows for differing amounts of sparsity of the vector $c$ in different levels.  Note that the levels $\mathbf{M}$ do not necessarily correspond to wavelet scales -- for now, we consider a general setting.  We also need a notion of best approximation:

\defn{[$(\mathbf{k},\mathbf{M}$)-term approximation]
Let $c$ be an element of either $\bbC^N$ or $\ell_2(\bbN)$.  We say that $c$ is 
$(\mathbf{k},\mathbf{M})$-compressible if $\sigma_{\mathbf{k},\mathbf{M}}(c)$ is small, where
\be{
\label{sigma_s_m}
\sigma_{\mathbf{k},\mathbf{M}}(c)_1 = \inf_{z \in \Sigma_{\mathbf{k},\mathbf{M}} } \| c - z \|_{1}.
}
}

As we have already seen for wavelet coefficients, it is often the case that $k_j / (M_j - M_{j-1}) \rightarrow 0$ as $j \rightarrow \infty$.  In this case, we say that $c$ is \textit{asymptotically sparse in levels}.  However, we stress that this framework does not explicitly require such decay.

We now consider the level-based sampling strategy:

\defn{[Multilevel random sampling]
\label{multi_level_dfn}
Let $r \in \bbN$, $\mathbf{N} = (N_1,\ldots,N_r) \in \bbN^r$ with $1 \leq N_1 
< \ldots < N_r$, $\mathbf{m} = (m_1,\ldots,m_r) \in \bbN^r$, with $m_j \leq 
N_j-N_{j-1}$, $j=1,\ldots,r$, and suppose that
\bes{
\Omega_j \subseteq \{ N_{j-1}+1,\ldots,N_{j} \},\quad | \Omega_j | = m_j,\quad 
j=1,\ldots,r,
}
are chosen uniformly at random, where $N_0 = 0$.  We refer to the set
\bes{
\Omega = \Omega_{\mathbf{N},\mathbf{m}} = \Omega_1 \cup \ldots \cup \Omega_r.
}
as an $(\mathbf{N},\mathbf{m})$-multilevel sampling scheme.
}

As discussed in Section \ref{ss:inc_CS}, the (infinite) Fourier/wavelets matrix $U = F \Phi$ is globally coherent.  However, as shown in Lemma \ref{l:DFT_Haar_local_coherence}, the coherence of its $(j,l)^{\rth}$ block is much smaller.  We therefore require a notion of local coherence:

\defn{[Local coherence]\label{loc_coherence}
Let $U$ be an isometry of either $\bbC^{N}$ or $\ell_2(\bbN)$.
If $\mathbf{N} = (N_1,\ldots,N_r) \in \bbN^r$ and $\mathbf{M} = 
(M_1,\ldots,M_r) \in \bbN^r$ with $1 \leq N_1 < \ldots N_r $ and $1 \leq M_1 < 
\ldots < M_r $ the $(j,l)^{\rth}$ local coherence of $U$ with respect to 
$\mathbf{N}$ and $\mathbf{M}$ is given by
\eas{
\mu_{\mathbf{N},\mathbf{M}}(j,l) &= 
\sqrt{\mu(P^{N_{j-1}}_{N_{j}}UP^{M_{l-1}}_{M_{l}})
\mu(P^{N_{j-1}}_{N_{j}}U)},\quad k,l=1,\ldots,r,
}
where $N_0 = M_0 = 0$ and $P^{a}_{b}$ denotes the projection matrix 
corresponding to indices $\{a+1,\hdots, b\}$.  In the case where $U$ is an operator on
$\ell_2(\bbN)$, we also define
\eas{
\mu_{\mathbf{N},\mathbf{M}}(j,\infty) &= 
\sqrt{\mu(P^{N_{j-1}}_{N_{j}}UP_{M_{r-1}}^\perp)
\mu(P^{N_{j-1}}_{N_{j}}U)},\quad j=1,\ldots,r.
}}

Note that the local coherence $\mu_{\mathbf{N},\mathbf{M}}(j,l)$ is not just the coherence $\mu(P^{N_{j-1}}_{N_{j}}UP^{M_{l-1}}_{M_{l}})$ in the $(j,l)^{\rth}$ block.  For technical reasons, one requires the product of this and the coherence $\mu(P^{N_{j-1}}_{N_{j}}U)$ in the whole $j^{\rth}$ row block.  

\begin{remark}
In \cite{KrahmerWardCSImaging}, the authors define a local coherence of a 
matrix $U$ in the $j^{\rth}$ row (as opposed to row block) to be the maximum 
of its entries in that row.  Using this, they prove recovery guarantees for 
the Fourier/Haar wavelets problem based on the RIP and the global sparsity 
$k$.  Unfortunately, these results do not explain the importance of structured 
sparsity as shown by the flip test.  Conversely, our notion of local coherence 
also takes into account the sparsity levels.  As we will see in Section 
\ref{ss:main_thms}, this allows one to establish recovery guarantees that 
are consistent with the flip test and properly explain 
the role of structure in the reconstruction.
\end{remark}

Recall that in practice (see Figure \ref{fig:U_plot}), the local coherence often decays along the diagonal blocks and in the off-diagonal blocks.  Loosely speaking, we say that the matrix $U$ is \textit{asymptotically incoherent} in this case.

In Section \ref{ss:Fourier_localinc} we argued that the Fourier/wavelets matrix was nearly block-diagonal.  In our theorems, we need to account for the off-diagonal terms.  To do this in the general setup, we require a notion of a relative sparsity:

\defn{[Relative sparsity]
\label{S}
Let $U$ be an isometry of either $\mathbb{C}^{N}$ or $\ell_2(\bbN)$.  For 
$\mathbf{N} = (N_1,\ldots,N_r) \in \bbN^r$, $\mathbf{M} = (M_1,\ldots,M_r) \in 
\bbN^r$ with $1 \leq N_1 < \ldots < N_r$ and $1 \leq M_1 < \ldots < M_r$, 
$\mathbf{k} = (k_1,\ldots,k_r) \in \bbN^r$ and $j=1,\ldots,r$, the 
$j^{\rth}$ relative sparsity is given by
\bes{
K_j = K_j(\mathbf{N},\mathbf{M},\mathbf{k}) =  \max_{\substack{z \in \Sigma_{\mathbf{k},\mathbf{M}}}, \| \eta \|_{\infty} \leq 1}\|P_{N_j}^{N_{j-1}}Uz\|^2,
}
where $N_0 = M_0 = 0$.
}

The relative sparsities $K_j$ take into account \textit{interferences}  between different sparsity level caused by the non-block diagonality of $U$.

\subsection{Main theorem}\label{ss:main_thms}

Given the matrix/operator $U$ and a multilevel sampling scheme $\Omega$, we now consider the solution of the convex optimization problem

\be{
\label{general_CS_U}
\min_{z \in \bbC^n} \| z \|_1\quad \mbox{s.t.}\quad \| P_{\Omega}y - P_{\Omega}U z \|_2 \leq \eta,
}
where $y = Uc + e$, $\| e \|_2 \leq \eta$.
Note that if $U = \Psi^* \Phi$, $x = \Phi c$ is the signal we wish to recover and $\hat{c}$ is a minimizer of \R{general_CS_U} then this gives the approximation $\hat{x} = \Phi \hat{c}$ to $x$.

\thm{
\label{main_full_fin_noise2}
Let $U \in \mathbb{C}^{N \times N}$ be an isometry and $c \in 
\mathbb{C}^{N}$.  Suppose that $\Omega = \Omega_{\mathbf{N},\mathbf{m}}$ is a 
multilevel sampling scheme, where $\mathbf{N} = (N_1,\ldots,N_r) \in \bbN^r$, 
$N_r = n$, and $\mathbf{m} = (m_1,\ldots,m_r) \in \bbN^r$.  Let $\epsilon \in (0,\E^{-1}]$ and suppose that 
$(\mathbf{k},\mathbf{M})$, where $\mathbf{M} = (M_1,\ldots,M_r) \in \bbN^r$,  
$M_r = n$, and $\mathbf{k} = (k_1,\ldots,k_r) \in \bbN^r$, are any pair such 
that the following holds: 
\begin{enumerate}
\item[(i)] We have
\be{
\label{bound_1}
1 \gtrsim \frac{N_j-N_{j-1}}{m_j}\log(\epsilon^{-1}) \left(
\sum_{l=1}^r \mu_{\mathbf{N},\mathbf{M}}(j,l) k_l\right)\log\left(n\right),\quad j=1,\ldots,r.
 }
 \item[(ii)] We have $m_j \gtrsim \hat m_j  \log(\epsilon^{-1}) \log\left(n\right)$, where $\hat{m}_j$ is such that
\be{
\label{bound_2}
 1 \gtrsim \sum_{j=1}^r \left(\frac{N_j-N_{j-1}}{\hat m_j} - 1\right) 
 \mu_{\mathbf{N},\mathbf{M}}(j,l)\tilde k_j,\quad l=1,\ldots,r,
 }
for all $\tilde{k}_1,\ldots,\tilde{k}_r \in (0,\infty)$ satisfying
\bes{
\tilde k_{1}+ \hdots + \tilde k_{r}  \leq k_1+ \hdots + k_r, \qquad \tilde k_j 
\leq K_j(\mathbf{N},\mathbf{M},\mathbf{k}).
}
\end{enumerate}
Suppose that $\hat{c} \in \bbC^N$ is a minimizer of \R{general_CS_U}.  Then, with probability exceeding $1-k\epsilon$,  where $k = k_1+\ldots+k_r$, 
we have that 
\be{
\label{eq:error}
\|c - \hat{c}\|_2\leq  C \left(\eta \sqrt{D} \left(1 +E\sqrt{k}\right) +\sigma_{\mathbf{k},\mathbf{M}}(c)_1 \right),}
for some constant $C$, where $\sigma_{\mathbf{k},\mathbf{M}}(c)_1$ is as in 
\R{sigma_s_m}, $
D= 1+ \frac{\sqrt{\log_2\left(6\epsilon^{-1}\right)}}{\log_2(4En\sqrt{s})}$ 
and $E = \max_{j=1,\ldots,r} \{ (N_{j}-N_{j-1})/m_j \}$.  If $m_j = 
N_{j}-N_{j-1}$, $j=1,\ldots,r$, then this holds with probability $1$.
}

A similar theorem can be stated and proved in the infinite-dimensional setting \cite{AHPRBreaking}.  For brevity, we shall not do this.

The key feature of Theorem \ref{main_full_fin_noise2} is that the bounds 
\R{bound_1} and \R{bound_2} involve only local quantities: namely, local 
sparsities $k_j$, local coherences $\mu(j,l)$, local measurements $m_j$ and 
relative sparsities $K_j$.  Note that the estimate \R{eq:error} directly 
generalizes a standard compressed sensing estimate for sampling with 
incoherent bases \cite{BAACHGSCS,Candes_Plan} to the case of multiple levels.  
Having said this, it is not immediately obvious how to understand these bounds 
in terms of how many measurements $m_j$ are actually required in the 
$j^{\rth}$ level.  However, it can be shown that these bounds are in fact 
sharp for a large class of matrices $U$ \cite{AHPRBreaking}.  Thus, little 
improvement is possible.  Moreover, in the important case Fourier/wavelets, 
one can analyze the local coherences $\mu(j,l)$ and relative sparsities $K_j$ 
to get such explicit estimates.  We consider this next.

\subsection{The case of Fourier sampling with wavelets}\label{ss:Fourier_wavelet}

Let us consider the example of Section \ref{ss:Fourier_localinc}, where the matrix $U$ arises from the Fourier/Haar wavelet pair, the sampling levels are correspond to the aforementioned frequency bands $W_j$ and the sparsity levels are the Haar wavelet scales.

\thm{\footnote{For a proof, we refer to \cite{AHRHaarFourierNote}.}
\label{t:Fourier_Haar}
Let $U$ and $\Omega$ be as in Section \ref{ss:Fourier_localinc} (recall that we index over $j,l=0,\ldots,r-1$) and suppose that $x \in \bbC^n$.  Let $\epsilon \in (0,\E^{-1}]$ and suppose that
\be{
\label{m_j_est}
m_j \gtrsim \left ( k_j + \sum^{r-1}_{\substack{l=0 \\ l \neq j}} 2^{-\frac{|j-l|}{2}} k_l \right )\log(\epsilon^{-1}) \log^2(n),\quad j=0,\ldots,r-1.
}
Then, with probability exceeding $1-k \epsilon$, where $k=k_0+\ldots+k_{r-1}$, any minimizer $\hat{x}$ of \R{general_CS} satisfies
\bes{
\| x - \hat{x} \|_2 \leq C \left ( \eta \sqrt{D} ( 1+E \sqrt{k} ) + \sigma_{\mathbf{k},\mathbf{M}}(\Phi^* x)_1 \right ),
}
where $\sigma_{\mathbf{k},\mathbf{M}}(\Phi^* x)_1$ is as in 
\R{sigma_s_m}, $
D= 1+ \frac{\sqrt{\log_2\left(6\epsilon^{-1}\right)}}{\log_2(4En\sqrt{k})}$ 
and $E = \max_{j=0,\ldots,r-1} \{ (N_{j}-N_{j-1})/m_j \}$.  If $m_j = |W_j|$, $j=0,\ldots,r-1$, then this holds with probability $1$.
}

The key part of this theorem is \R{m_j_est}.  Recall that if $U$ were exactly block diagonal, then $m_j \gtrsim k_j$ would suffice (up to log factors).  The estimate \R{m_j_est} asserts that we require only slightly more samples, and this is due to interferences from the other sparsity levels.  However, as $|j-l|$ increases, the effect of these levels decreases exponentially.  Thus, the number of measurements $m_j$ required in the $j^{\rth}$ frequency band is determined predominantly by the sparsities in the scales $l \approx j$.  Note that $k_l \approx \ord{k_j}$ when $l \approx j$ for typical signals and images, so the estimate \R{m_j_est} is typically on the order of $k_j$ in practice.

The estimate \R{m_j_est} both agrees with the conclusion of the flip test in 
Figure \ref{f:BernStoneHad} and explains the results seen.  Flipping the 
wavelet coefficients changes the local sparsities $k_1,\ldots,k_r$.  Therefore 
to recover the flipped image to the same accuracy as the unflipped image, 
\R{m_j_est} asserts that one must change the local numbers of measurements 
$m_j$.  But in Figure \ref{f:BernStoneHad} the same sampling pattern was used 
in both cases, thereby leading to the worse reconstruction in the flipped 
case.  Note that \R{m_j_est} also highlights why the optimal sampling pattern 
must depend on the image, and specifically, the local sparsities.  In 
particular, there can be no optimal sampling strategy for all images. 

Note that Theorem \ref{t:Fourier_Haar} is a simplified version, presented here for the purposes of elucidation, of a more general result found in \cite{AHPRBreaking} which applies to all compactly-supported orthonormal wavelets in the infinite-dimensional setting.

\section{Structured sampling and structured recovery}\label{ss:Structured_Sampling_Recovery}


\begin{figure}
\centering
\begin{tabular}{@{}c@{\hspace{0.02\textwidth}}c@{}}
\includegraphics[width=0.49\textwidth]{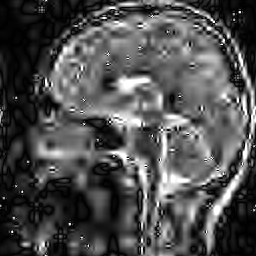}&
\includegraphics[width=0.49\textwidth]{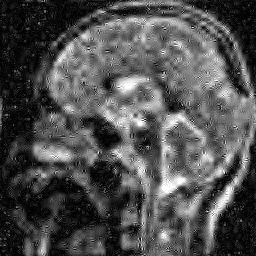}\\
Rand. Gauss to DB4, ModelCS, err = 41.8\% & Rand. Gauss to DB4, TurboAMP, err 
= 39.3\% \\[7pt]
\includegraphics[width=0.49\textwidth]{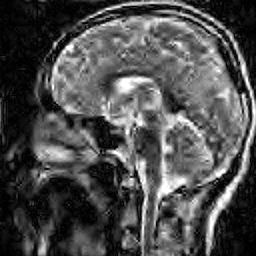}&
\includegraphics[width=0.49\textwidth]{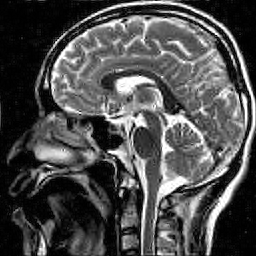}\\
Rand. Gauss to DB4, BCS, err = 29.6\% & Multilevel DFT to DB4, $\ell^1$, err = 
18.2\%
\end{tabular}
\caption{Recovery from 12.5\% measurements at \res{256}. Comparison between 
random sampling with structured recovery and structured sampling with 
$\ell^1$-minimization recovery.}
\label{fig:Comp4SamplVsRec}       
\end{figure}

Structured sparsity within the context of compressed sensing has been 
considered in numerous previous works.  See 
\cite{BaranuikModelCS,BourrierEtAlInstance,CalderbankCommInspCS,donohoCS,DuarteEldarStructuredCS,EldarUnionSubspace,HeCarinStructCS,HeCarinTreeCS,EldarXampling2,DonohoCSext,CalderbankNIPS,TurboAMP}
 and references therein.  For the problem of reconstructing wavelet 
coefficients, most efforts have focused on their inherent tree structure (see 
Remark \ref{r:tree_structure}).  Three well-known algorithmic approaches for 
doing this are model-based compressed sensing \cite{BaranuikModelCS}, TurboAMP \cite{TurboAMP}, and Bayesian compressed sensing \cite{HeCarinStructCS,HeCarinTreeCS}. 
All methods use Gaussian or Bernoulli random measurements, and seek to leverage the wavelet tree structure -- the former deterministically, the latter two in a probabilistic manner -- by appropriately-designed recovery algorithms (based on modifications of existing iterative algorithms for compressed sensing).  In other words, structure is incorporated solely in the recovery algorithm, and not in the measurements themselves.

In Figure \ref{fig:Comp4SamplVsRec} we compare these algorithms with the 
previously-described method of multilevel Fourier sampling (similar results are also witnessed with the Hadamard matrix). Note 
that the latter, unlike other three methods, exploits structure by taking 
appropriate measurements, and uses an unmodified compressed sensing algorithm 
($\ell^1$ minimization).  As is evident, this approach is able to better 
exploit the sparsity structure, leading to a significantly improved 
reconstruction.  This experiment is representative of a large set tested. In 
all cases, we find that exploiting structure by sampling with asymptotically 
incoherent Fourier/Hadamard bases outperforms such approaches that seek to 
leverage structure in the recovery algorithm.

\section{The Restricted Isometry Property in levels}
\begin{figure}[t]
\small
\begin{center}
\begin{tabular}{@{}c@{\hspace{0.005\linewidth}}c@{\hspace{0.005\linewidth}}c@{}}
\includegraphics[width=0.33\linewidth]{0256_12p5_dft_rwtdwtDB4_5_VQP9_map.png}&
\includegraphics[width=0.33\linewidth]{0256_12p5_dft_rwtdwtDB4_5_VQP9_rec.png}&
\includegraphics[width=0.33\linewidth]{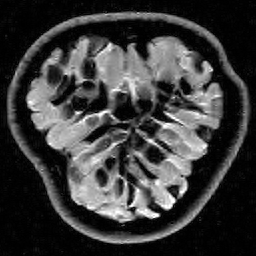}\\
Subsampling pattern & DFT to DB4, & DFT to Levels Flipped DB4, \\
  & error = 10.95\% & error = 12.28\%\\
\end{tabular}
\caption{The flip test in levels for the image considered in Figure \ref{f:BernStoneHad}.}
\label{f:BernStoneHad2}
\end{center}
\end{figure}

The flip test demonstrates that the subsampled Fourier/wavelets matrix $P_{\Omega}U$ does not satisfy a meaningful RIP.  However, due to the sparsity structure of the signal, we are still able to reconstruct, as was confirmed in Theorem \ref{main_full_fin_noise2}.  The RIP is therefore too crude a tool to explain the recoverability properties of structured sensing matrices.   Having said that, the RIP is a useful for deriving uniform, as opposed to nonuniform, recovery results, and for analyzing other compressed sensing algorithms besides convex optimization, such as greedy methods \cite{FoucartRauhutCSbook}.  This raises the question of whether there are alternatives which are satisfied by such matrices.  One possibility which we now discuss is the RIP in levels.  Throughout this section we shall work in the finite-dimensional setting.  For proofs of the theorems, see \cite{Bastounis}.

\begin{definition}\label{def:RIPInLevels}
Given an $r$-level sparsity pattern $(\mathbf{k},\mathbf{M})$, where $M_r =n$, we say that the matrix $A \in \bbC^{m \times n}$ satisfies the \emph{RIP in levels} ($\text{RIP}_L$) with $\text{RIP}_L$ constant $\delta_{\mathbf{k}} \geq 0$ if for all $x$ in $\Sigma_{\mathbf{k,M}}$ we have
\begin{equation*}
(1-\delta_{\mathbf{k}})\| x \|^2_2 \leq \| A x \|^2_2 \leq (1+\delta_{\mathbf{k}}) \| x \|^2_2.
\end{equation*}
\end{definition}
The motivation for this definition is the following.  Suppose that we repeat the flip test from Figure \ref{f:BernStoneHad} except that instead of completely flipping the coefficients we only flip them within levels corresponding to the wavelet scales. We will refer to this as the \textit{flip test in levels}. Note the difference between Figure  \ref{f:BernStoneHad} and Figure  \ref{f:BernStoneHad2}, where latter presents the flip test in levels: clearly, flipping within scales does not alter the reconstruction quality.  In light of this experiment, we propose the above RIP in levels definition so as to respect the level structure.

We now consider the recovery properties of matrices satisfying the RIP in levels.  For this, we define the \textit{ratio constant} $\lambda_{\mathbf{k,M}} $ of a sparsity pattern $(\mathbf{k},\mathbf{M})$ to be 
$
\lambda_{\mathbf{k,M}} := \max_{j,l} k_j / k_l.
$
We assume throughout that $k_j \geq 1$, $\forall j$, so that $\eta_{\mathbf{k},\mathbf{M}} < \infty$, and also that $M_r = n$. 
We now have the following:

\begin{theorem}\label{Theorem:RIPInLevelsRecoveryTheorem}
Let $(\mathbf{k},\mathbf{M})$ be a sparsity pattern with $r$ levels and ratio constant $\lambda = \lambda_{\mathbf{k,M}}$. Suppose that the matrix $A$ has $\text{RIP}_L$ constant $\delta_{2\mathbf{k}}$ satisfying
\begin{equation}\label{delta_bound}
\delta_{2\mathbf{k}} < \frac{1}{\sqrt{r (\sqrt{\lambda}+1/4)^2 + 1}}.
\end{equation}
Let $x \in \bbC^n$, $y \in \bbC^m$ be such that $\| U x - y \|_2 \leq \eta$, and let $\hat{x}$ be a minimizer of 
\bes{
\min_{z \in \bbC^n} \| z \|_1\quad \mbox{s.t.}\quad \| y - A z \|_2 \leq \eta.
}
Then 
\begin{equation}\label{recovery_bound} 
\| x - \hat{x} \|_1 \leq C \sigma_{\mathbf{k},\mathbf{M}}(x)_1 + D \sqrt{k}\eta ,
\end{equation}
where $k =  k_1+\ldots+k_r$ and the constants $C$ and $D$ depend only on $\delta_{2\mathbf{k}}$. 
\end{theorem}

This theorem is a generalization of a known result in standard (i.e.\ one-level) compressed sensing.  Note that \R{delta_bound} reduces to the well-known estimate $\delta_{2 k} \leq 4 / \sqrt{41}$ \cite{FoucartRauhutCSbook} when $r=1$.  On the other hand,  in the multiple level case the reader may be concerned that the bound ceases to be useful, since the right-hand side of \R{delta_bound} deteriorates with both the number of levels $r$ and the sparsity ratio $\lambda$.  As we show in the following two theorems, the dependence on $r$ and $\lambda$ in \R{delta_bound} is sharp:

\begin{theorem}\label{Theorem:etaDependenceTheorem}
Fix $a \in \mathbb{N}$. There exists a matrix $A$ with two levels and a sparsity pattern $(\mathbf{k},\mathbf{M})$ such that the $\text{RIP}_L$ constant $\delta_{a\mathbf{k}}$ and ratio constant $\lambda = \lambda_{\mathbf{k,M}}$ satisfy
\begin{equation}\label{equation:etaDependenceRIPBounds}
\delta_{a\mathbf{k}} \leq 1/|f(\lambda)|,
\end{equation}
where $f(\lambda) = o(\sqrt{\lambda})$, but there is an $x \in \Sigma_{\mathbf{k},\mathbf{M}}$ such that $x$ is not the minimizer of 
\bes{
\min_{z \in \bbC^n} \| z \|_1\quad \mbox{s.t.}\quad A z = A x .
}
\end{theorem}
Roughly speaking, Theorem \ref{Theorem:etaDependenceTheorem} says that if we fix the number of levels and try to replace \R{delta_bound}
with a condition of the form
\begin{equation*}
\delta_{2\mathbf{k}} < \frac{1}{C\sqrt{r}}\lambda^{-\frac{\alpha}{2}}
\end{equation*}
for some constant $C$ and some $\alpha < 1$ then the conclusion of Theorem \ref{Theorem:RIPInLevelsRecoveryTheorem} ceases to hold. In particular, the requirement on $\delta_{2\mathbf{k}}$ cannot be independent of $\lambda$. The parameter $a$ in the statement of Theorem \ref{Theorem:etaDependenceTheorem} also means that we cannot simply fix the issue by changing $\delta_{2\mathbf{k}}$ to $\delta_{3\mathbf{k}}$, or any further multiple of $\mathbf{k}$.

Similarly, we also have a theorem that shows that the dependence on the number of levels $r$ cannot be ignored.
\begin{theorem}\label{Theorem:lDependenceTheorem}
Fix $a \in \mathbb{N}$. There exists a matrix $A$ and an $r$-level sparsity pattern $(\mathbf{k},\mathbf{M})$ with ratio constant $\lambda_{\mathbf{k,M}} = 1$ such that the $\text{RIP}_L$ constant $\delta_{a\mathbf{k}}$ satisfies
\begin{equation*}
\delta_{a\mathbf{k}} \leq 1/|f(r)|,
\end{equation*}
where $f(r) = o(\sqrt{r})$, but there is an $x \in \Sigma_{\mathbf{k},\mathbf{M}}$ such that $x$ is not the minimizer of 
\bes{
\min_{z \in \bbC^n} \| z \|_1\quad \mbox{s.t.}\quad A z = A x .
}
\end{theorem}

These two theorems suggests that, at the level of generality of the $\text{RIP}_L$, one must accept a bound that deteriorates with the number of levels and ratio constant.  This begs the question: what is the effect of such deterioration?  To understand this, consider the case studied earlier, where the $r$ levels correspond to wavelet scales.  For typical images, the ratio constant $\lambda$ grows only very mildly with $n$, where $n = 2^r$ is the dimension.  Conversely, the number of levels is equal to $\log_2(n)$.  This suggests that estimates for the Fourier/wavelet matrix that ensure an RIP in levels (thus guaranteeing uniform recovery) will involve at least several additional factors of $\log(n)$ beyond what is sufficient for nonuniform recovery (see Theorem \ref{t:Fourier_Haar}).  Proving such estimates for the Fourier/wavelets matrix is work in progress.

\begin{acknowledgement}
The authors thank Andy Ellison from Boston University Medical School for 
kindly providing the MRI fruit image, and General Electric Healthcare for 
kindly providing the brain MRI image.
BA acknowledges support from the NSF DMS grant 1318894. ACH acknowledges 
support from a Royal Society University Research Fellowship. ACH and BR 
acknowledge the UK 
Engineering and Physical Sciences Research Council (EPSRC) grant EP/L003457/1.
\end{acknowledgement}

\bibliographystyle{plain}
\bibliography{BookChptRefs}
\end{document}